\documentclass[reqno,12pt]{amsart}

\usepackage[T1]{fontenc}
\usepackage{lmodern}

\usepackage[numbers,sort&compress]{natbib}
\usepackage{graphicx}
\usepackage{amsmath,amssymb}
\usepackage{amsaddr}
\usepackage{mathtools}
\usepackage{textcomp}
\usepackage{algorithm}
\usepackage{algpseudocode}
\usepackage[labelfont={rm}, justification=centering]{subfig}
\usepackage{booktabs}

\calclayout

\renewcommand{\vec}[1]{\boldsymbol{#1}}

\newcommand{\D}    {\mathrm{d}}

\newcommand{\NM}[1]{\underline{#1}}

\newcommand{\jump}[1]%
  {[\hspace{-0.15em}[\hspace{0.05em}#1\hspace{0.05em}]\hspace{-0.15em}]}

\newcommand{\average}[1]%
  {\{\hspace{-0.20em}\{\hspace{0.05em}#1\hspace{0.05em}\}\hspace{-0.20em}\}}

\newcommand{\transpose}[1]{#1^{\textsc{t}}}

\newcommand{\N}      [1]   {N_{\textsc{#1}}}
\newcommand{\ns}     [1][] {N_{\textsc{s}#1}}
\newcommand{\deltaO} [1][] {\delta_{\textsc{o}#1}}
\newcommand{\no}     [1][] {N_{\textsc{o}#1}}

\newcommand{\IOP}{\NM{\mathcal I}}

\begin{document}

\title%
[% Running title
Robust MG-DG for the Poisson Equation in 3D
]
{% Full title
Robust Multigrid for Cartesian Interior Penalty DG Formulations of the Poisson Equation in 3D
}
\author{J\"org Stiller}
\address{%
  Technische Universit\"at Dresden, Institute of Fluid Mechanics and
  Center for Advancing Electronics Dresden, 01062 Dresden, Germany}
\email{joerg.stiller@tu-dresden.de}

\begin{abstract}
We present a polynomial multigrid method for the nodal interior penalty formulation
of the Poisson equation on three-dimensional Cartesian grids.
Its key ingredient is a weighted overlapping Schwarz smoother operating on
element-centered subdomains.
The MG method reaches superior convergence rates
corresponding to residual reductions of about two orders of magnitude within a single
V(1,1) cycle.
It is robust with respect to the mesh size and the ansatz order, at least up to
${P=32}$.
Rigorous exploitation of tensor-product factorization yields a computational
complexity of $O(PN)$ for $N$ unknowns, whereas numerical experiments indicate even
 linear runtime scaling.
Moreover, by allowing adjustable subdomain overlaps and adding Krylov acceleration,
the method proved feasible for anisotropic grids with element aspect ratios up to 48.
\end{abstract}


\maketitle

%===============================================================================

\section{Introduction}
\label{sec:intro}

Discontinuous Galerkin (DG) methods combine multiple desirable properties of
finite element and finite volume methods, including geometric flexibility,
variable approximation order, straightforward adaptivity and suitability for
conservation laws \cite{Cockburn2000,Hesthaven_NDGM_2008}.
Though traditionally focused on hyperbolic systems, the need for implicit diffusion schemes and application to other problem classes, such as incompressible flow and elasticity, led to a growing interest in DG methods and related solution techniques
for elliptic equations \cite{Arnold2001,Riviere2008}.
This paper is concerned with fast elliptic solvers based on the multigrid (MG)
method.
In the context of high-order spectral element and DG methods, several approaches
have been proposed:
polynomial or $p$-MG \cite{HAM03,FOLD05,HA06,HSN13},
geometric or $h$-MG \cite{GK03,Kan04,KT08},
and
algebraic MG \cite{OS11,BBS12}.
The most efficient methods reported so far \cite{BBS12,Stiller2016}
use block smoothers that can be regarded as overlapping Schwarz methods.
This work presents a hybrid Schwarz/MG method for nodal interior penalty DG
formulations of Poisson problems on 3D Cartesian grids.
It extends the techniques put forward in \cite{Stiller2016,Stiller2017}
and generalizes the approach to variable subdomain overlaps.
The remainder of the paper is organized as follows:
Section~\ref{sec:dg} briefly describes the discretization,
Sec.~\ref{sec:multigrid} the multigrid technique, including the Schwarz smoother,
and Sec.~\ref{sec:results} presents the results of the assessment by means of
numerical experiments.
Section~\ref{sec:concl} concludes the paper.

%===============================================================================

\section{Discontinuous Galerkin Method}
\label{sec:dg}

We consider the Poisson equation
\begin{equation}
\label{eq:poisson}
-\nabla^2 u = f
\end{equation}
in the periodic domain
${\Omega  = [0,l_x] \times [0,l_y] \times [0,l_z]}$.%
\footnote{Within this paper, the following symbols are used concurrently for representing the Cartesian coordinates:
$\vec{x} = [x_i] = (x_1, x_2, x_3) = (x,y,z)$.}
For discretization, the domain is decomposed into cuboidal elements $\{\Omega^e\}$ forming a Cartesian mesh.
The discrete solution $u_h$ is sought in the function space
\begin{equation}
\label{eq:function space}
\mathcal V_h =
  \left\{
    v \in L^2(\Omega) : v|_{\Omega^{e}} \in \mathcal P_P^3(\Omega^e) \quad
    \forall \Omega^e \subset \Omega
  \right\},
\end{equation}
where $\mathcal P_P^3$ is the 3D tensor product of polynomials of at most degree $P$.
To cope with discontinuity we introduce the interior surface ${\Gamma = \cup \Gamma^f}$
which is composed of the element interfaces $\{\Gamma^f\}$.
Let $\Omega^-$ and $\Omega^+$ denote the elements adjacent to $\Gamma^f$ and, respectively,
$\vec{n}^-$ and $\vec{n}^+$ their exterior normals, and
$v^-$ and $v^+$ the restrictions of $v$ to the joint face from inside the elements.
Then we define the average and jump operators as
\begin{equation}
\average{v} = \frac{1}{2}(v^- + v^+),
\quad
\jump{v} = \vec{n}^- v^- + \vec{n}^+v^+
.
\end{equation}
Given these prerequisites, the interior penalty discontinuous Galerkin formulation can be stated as follows (see, e.g. \cite{Arnold2001}):
Find ${u_h \in \mathcal V_h}$ such that for all ${v \in \mathcal V_h}$

\begin{multline}
\label{eq:ip}
\int_{\Omega} \nabla v \cdot \nabla u_h \D\Omega
+ \int_{\Gamma}
    \left( \jump{\nabla v}\cdot\average{u_h}
        + \average{v}\cdot\jump{\nabla u_h}
    \right) \D\Gamma
\\
+ \int_{\Gamma} \mu \jump{v} \cdot \jump{u_h}
= \int_{\Omega} v f \D\Omega
\, ,
\end{multline}
where $\mu$ is a piece-wise constant penalty parameter that must be chosen large enough to ensure stability.

Although, in theory, any suitable basis in $\mathcal V_h$ can be chosen, we restrict ourselves  to nodal tensor-product bases generated from the Lagrange polynomials to the Gauss-Legendre (GL) or Gauss-Legendre-Lobatto (GLL) points, respectively.
The discrete solution is expressed in $\Omega^e$ as
\begin{equation}
\label{eq:u_h}
u_h(\vec{x})|_{\Omega^{e}}
= u^e(\vec{\xi}^e(\vec{x}))
= \sum_{i,j,k=0}^P u_{ijk}^e \varphi_i(\xi) \varphi_j(\eta) \varphi_k(\zeta)
\, ,
\end{equation}
where
$\varphi_i$ denotes the 1D base functions and
$\vec{\xi}^e(\vec{x})$ the transformation from
${\Omega^e}$ to the reference element
${[-1,1]^3}$.
Each coefficient $u_{ijk}^e$ is associated with one local base function, which is globalized by zero continuation outside $\Omega^e$.
Inserting these base functions in \eqref{eq:ip} for $v$ and applying GL or GLL quadrature, according to the chosen basis, yields the discrete equations
\begin{equation}
\label{eq:discrete}
\NM A\, \NM u = \NM f
\end{equation}
for the solution vector ${\NM u = [u_{ijk}^e]}$.
Due to the Cartesian element mesh and the tensor-product ansatz \eqref{eq:u_h} the system matrix assumes the tensor-product form
\begin{equation}
\label{eq:A}
\NM A = \NM M_z \otimes \NM M_y \otimes \NM L_x
      + \NM M_z \otimes \NM L_y \otimes \NM M_x
      + \NM L_z \otimes \NM M_y \otimes \NM M_x
\,
\end{equation}
where
$\NM M_{\ast}$, $\NM L_{\ast}$ are the 1D mass and stiffness matrices for
${\ast=x,y,z}$.
Without going into detail we remark that
$\NM M_{\ast}$ is positive diagonal, and
$\NM L_{\ast}$ symmetric positive semi-definite and block tridiagonal
for either basis choice.
The rigorous exploitation of these properties is crucial for the efficiency of the overall method.

%===============================================================================

\section{Multigrid Techniques}
\label{sec:multigrid}

The tensor-product structure of \eqref{eq:discrete} allows for a straight-forward extension of the multigrid techniques developed in \cite{Stiller2016} for the 2D case.
In the following, we examine polynomial multigrid (MG) and multigrid-preconditioned conjugate gradients (MG-CG) both using an overlapping Schwarz method for smoothing.

%-------------------------------------------------------------------------------

\subsection{Schwarz Method}
\label{sec:multigrid:schwarz}

\begin{figure}%[b]
\includegraphics[width=60mm]{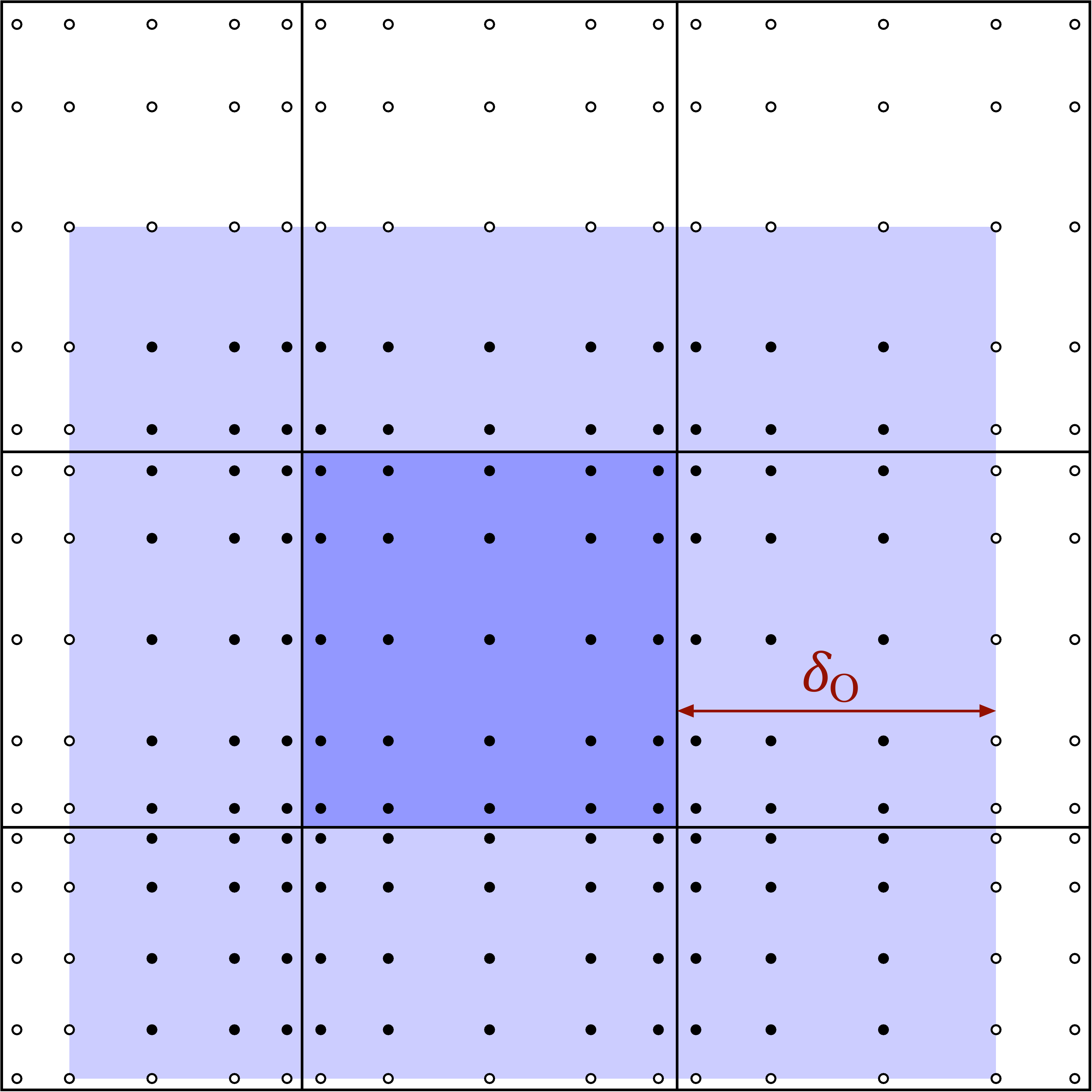}
\caption{Element-centered subdomain.
Every subdomain consists of a core region coinciding with the embedded element (dark) and an overlap zone (light shaded).
The latter represents a strip of variable width $\deltaO$, which is adopted from the surrounding elements.
The circles are the GL nodes for polynomial order ${P=4}$.
Filled circles indicates the unknowns that are solved for.}
\label{fig:subdomain}
\end{figure}

Schwarz methods are iterative domain decomposition techniques which improve the approximate solution by parallel or sequential subdomain solves, leading to additive or multiplicative methods, respectively.
Here, we consider additive Schwarz with overlapping element-centered subdomains as sketched in Fig.~\ref{fig:subdomain}.
The overlap width $\deltaO$ can be different on each side of the embedded element, but may not exceed the width of the adjoining element. Alternatively, the overlap can be specified by prescribing the number $\no$ of node layers adopted from the latter.

To derive the subdomain problems, we first rewrite \eqref{eq:discrete} into the residual form
\begin{equation}
\label{eq:discrete}
\NM A\, \Delta \NM u = \NM f - \NM A\, \NM{\tilde u} = \NM r
\,,
\end{equation}
where
${\Delta \NM u = \NM u - \NM{\tilde u}}$ is the correction to the current approximate solution
$\NM{\tilde u}$.
For each subdomain $\Omega_s$ we define the restriction operator $\NM R_s$ such that
${\NM u_s = \NM R_s \NM u}$ gives the associated coefficients.
Conversely, the transposed restriction operator, $\transpose{\NM R_s}$ globalizes the local coefficients by adding zeros for exterior nodes.
With these prerequisites the correction contributed by $\Omega_{s}$ is defined as the solution to the subproblem
\begin{equation}
\label{eq:subproblem}
\NM A_{ss} \Delta \NM u_s = \NM r_s \, ,
\end{equation}
where
\mbox{$\NM A_{ss} = \NM R_s \NM A \, \transpose{\NM R}_s$}
is the restricted system matrix and
\mbox{$\NM r_s = \NM R_s \NM{r}$}
the restricted residual.
Due to the cuboidal shape of the subdomain, the restriction operator possesses the tensor-product factorization
${\NM R_s = \NM R_{s,x} \otimes \NM R_{s,y} \otimes \NM R_{s,z}}$
and, as a consequence, $\NM A_{ss}$ inherits the tensor-structure of the full system matrix
\eqref{eq:A}.
Moreover it is regular and can be inverted using the fast diagonalization technique of Lynch et al. \cite{Lyn64} to obtain
\begin{equation*}
\NM A_{ss}^{-1}
= ( \NM S_z \otimes \NM S_y \otimes \NM S_x )
  ( \NM I \otimes \NM I \otimes \NM \Lambda_x
  + \NM I \otimes \NM \Lambda_y \otimes \NM I
  + \NM \Lambda_z \otimes \NM I \otimes \NM I )^{-1}
  ( \transpose{\NM S_z} \otimes \transpose{\NM S_y} \otimes \transpose{\NM S_x} ),
\end{equation*}
where
$\NM S_{\ast}$ is the column matrix of eigenvectors and
$\Lambda_{\ast}$ the diagonal matrix of eigenvalues
to the generalized eigenproblem for the restricted 1D stiffness and mass matrices and
${\ast=x,y,z}$.
Exploiting this structure the solution can be computed in $O(P^4)$ operations per subdomain.

\begin{figure}%[b]
\includegraphics[width=72mm]{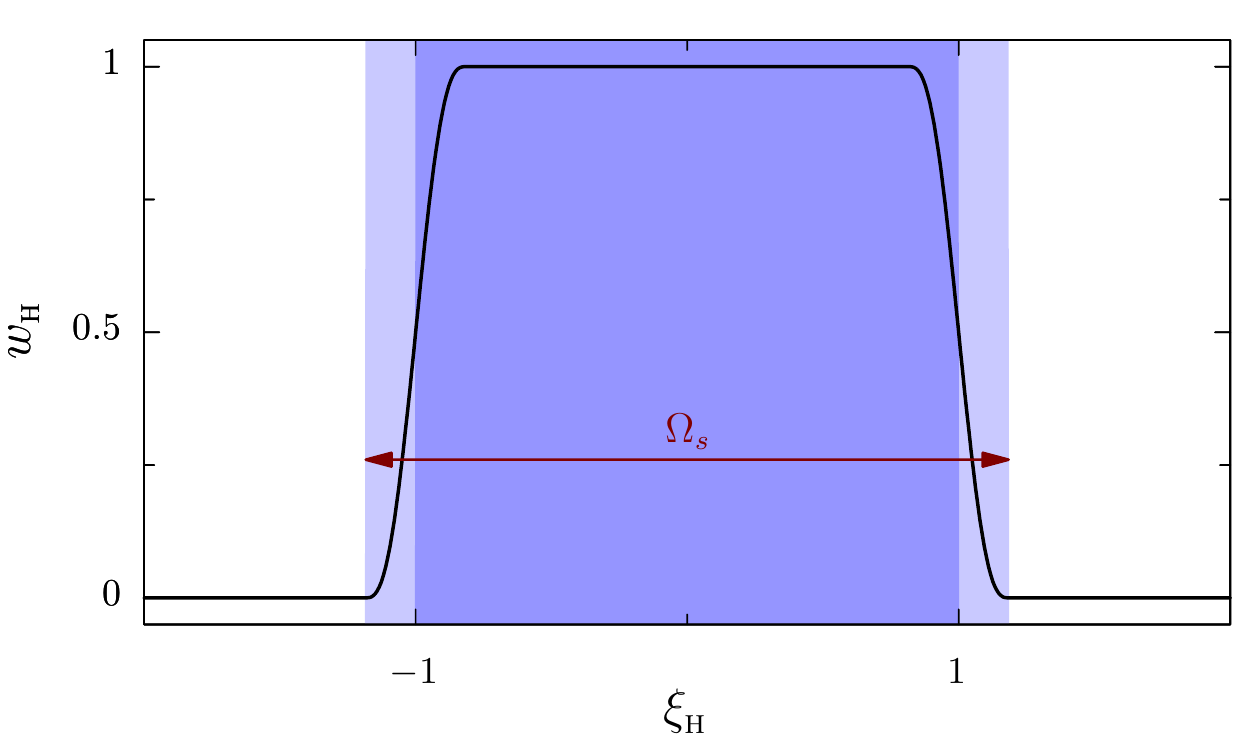}
\caption{Hat-shaped weight function using a piece-wise quintic transition from 0 at the subdomain boundary to 1 in the non-overlapped part of the core region.
The coordinate $\xi_{\textsc h}$ coincides with $\xi$ inside the central element, and $\xi\mp2$ in the left and right neighbor elements, respectively.
}
\label{fig:wH}
\end{figure}

One additive Schwarz iteration proceeds as follows:
First, all subproblems are solved in parallel, which yields the local corrections
${\Delta \NM u_s}$.
Afterwards, the global correction is computed as the weighted average
\begin{equation}
\label{eq:correction:global}
\Delta \NM u \simeq \sum_{s} \transpose{\NM R_{s}} (\NM W_s \Delta \NM u_{s}) \, ,
\end{equation}
where
${\NM W_s = \NM W_z \otimes \NM W_y \otimes \NM W_x}$
is the diagonal local weighting matrix.
The constituent 1D weights are computed from the hat-shaped weight function $w_{\textsc h}$, which is  illustrated in Fig.~\ref{fig:wH}.
The complete definition of the weight function and alternative choices are given in \cite{Stiller2016}.

%-------------------------------------------------------------------------------

\subsection{Multigrid and Preconditioned Conjugate Gradient Methods}
\label{sec:multigrid:pmg}

For MG we define a series of polynomial levels $\{P_l\}$
with
\mbox{$P_l = 2^l$}
increasing from $1$ at \mbox{$l=0$} to $P$ at top level $L$.
Correspondingly, let
$\NM u_l$ denote the global coefficients and
$\NM A_l$ the system matrix on level $l$.
On the top level we have
\mbox{$\NM u_L = \NM u$} and
\mbox{$\NM A_L = \NM A$},
whereas on lower levels $\NM u_l$ is the defect correction and $\NM A_l$ the counterpart of $\NM A$ obtained with elements of order $P_l$.
For transferring the correction from level \mbox{$l-1$} to level $l$ we use the embedded interpolation operator $\IOP_l$, and for restricting the residual its transpose.
These ingredients allow to build the multigrid V-cycle summarized in Algorithm~\ref{alg:V-cycle}, where the \textsc{Smoother} represents the weighted additive Schwarz method.
To allow for variable V-cycles \cite{Bra93}, the number of pre- and post-smoothing steps, $\ns[1,l]$ and $\ns[2,l]$, may change from level to level.
Line~\ref{alg:V-cycle:coarse} of Algorithm~\ref{alg:V-cycle} defines the coarse grid solution formally by means of the pseudoinverse $\NM A_0^{+}$.
In our implementation the coarse problem is solved using the conjugate gradient method. To ensure convergence, the right side is projected to the null space of $\NM A_0$, as proposed in \cite{Kaa88}.

\begin{algorithm}[ht]
\caption{Multigrid V-cycle.}
\label{alg:V-cycle}
\begin{algorithmic}[1]
\Function{MultigridCycle}{$\NM u$, $\NM f$, $\NM\ns$}
   \State $\NM u_L \gets \NM u$
   \State $\NM f_L \gets \NM f$
   \For{$l=L,1$ \textbf{step} $-1$}
       \If{$l < L$}
          \State $\NM u_l \gets 0$
       \EndIf
       \State $\NM u_l \gets$ \Call{Smoother}{$\NM u_l$, $\NM f_l$, $\ns[1,l]$}
         \Comment{Pre-smoothing}
       \State $\NM f_{l-1} \gets \transpose{\IOP}_l (\NM f_l - \NM A_l \NM u_l)$
         \Comment{Residual restriction}
   \EndFor
   \State $\NM u_0 \gets \NM A_0^{+} \NM f_0$
     \label{alg:V-cycle:coarse}
     \Comment{Coarse grid solution}
   \For{$l=1,L$}
       \State $\NM u_l \gets \NM u_l + \IOP_l \NM u_{l-1}$
         \Comment{Correction prolongation}
       \State $\NM u_l \gets$ \Call{Smoother}{$\NM u_l$, $\NM f_l$, $\ns[2,l]$}
         \Comment{Post-smoothing}
   \EndFor
   \State \textbf{return} $\NM u \gets \NM u_L$
\EndFunction
\end{algorithmic}
\end{algorithm}

It is well known that the robustness of multigrid method can be enhanced by Krylov acceleration \cite{TOS00}.
Here we use the inexact preconditioned conjugate gradient method \cite{GY99}, which copes with the asymmetry introduced by the weighted Schwarz method without imposing significant extra cost in comparison to conventional CG.
A detailed description of the algorithm is given in \cite{Stiller2016}.

%===============================================================================

\section{Results}
\label{sec:results}

For assessing robustness and efficiency, the described methods were implemented in Fortran and applied to
the $2\pi$-periodic Poisson problem with the exact solution
\[
  u(\vec{x}) = \cos(x - 3x + 2z)
      \sin(1 + x)
      \sin(1 - x)
      \sin(2x + x)
      \sin(3x - 2y + 2z).
\]
To keep the test series manageable, we constrained ourselves to equidistant grids with an identical number of elements in each direction.
Anisotropic meshes were realized be choosing the domain extensions as multiples of $2\pi$, i.e.
${l_{\ast} = 2 \pi s_{\ast}}$, which yields the aspect ratio
${\Delta x : \Delta y : \Delta z = s_x:s_y:s_z}$.
All tests started from a random guess confined to ${[-1,1]}$
and used a penalty parameter of ${\mu_{\ast}=2\mu_{\min,\ast}}$,
where $\mu_{\min,\ast}$ is the stability threshold,
e.g., $\mu_{\min,x} = P(P+1)/\Delta x$ for the $x$ direction
\cite{Stiller2016}.
The program was compiled using the Intel Fortran compiler 17.0 with optimization \textminus O3 and run on a 3.1\,GHz Intel Core i7-5557U CPU.

The primary assessment criterium is the average multigrid convergence rate
\begin{equation*}
\label{eq:rho}
\rho = \sqrt[n]{\frac{r_n}{r_0}}
\, ,
\end{equation*}
where $r_n$ is the Euclidean norm of the residual vector after the $n$th cycle.
Additionally we consider the number of cycles $n_{10}$ and the average runtime per unknown $\tau_{10}$ that are required to reduce the residual by a factor of $10^{10}$.
These quantities follow from the convergence rate by
${n_{10}} = \lceil -10/\lg \rho \rceil$ and
${\tau_{10} = -10\, t_{\textsc{c}}/\lg \rho}$,
respectively, where
$t_{\textsc{c}}$ is the time required for one V-cycle.

%-------------------------------------------------------------------------------

\subsection{Isotropic Meshes}

First we consider the isotropic case with ${s_x=s_y=s_z=1}$, such that
${\Omega = [0,2\pi]^3}$.
For assessing the impact of the subdomain overlap on the convergence rate and
computational cost, we performed a test series for ansatz orders ${P=4}$ to $32$
using a degree-dependent tessellation into ${\N{e}=(128/P)^3}$ elements.
Table~\ref{tab:isotropic:overlap} presents the logarithmic convergence rates
for 14 test cases featuring different choices for
the basis functions (GLL or GL),
the solution method (MG or MG-CG)
and the subdomain overlap.
Note that the latter was chosen identical in each direction, because of mesh isotropy.
All cases employed a fixed V-cycle with one pre- and post-smoothing step.
Independent of the basis and the solution method, choosing a minimal overlap
of just one node (${\no=1}$) yields acceptable convergences rates for low order (${P=4}$), but becomes inefficient with increasing order.
Using a geometrically fixed overlap of just eight percent of the element width (${\deltaO = 0.08\Delta x}$) on
every mesh level ensures robustness with respect to the ansatz order and even
improves the convergence with growing $P$.
Enlarging the overlap increases the convergence rate but also the computational cost,
as will be detailed in a moment.
Using Krylov acceleration tends to give faster convergence, however, this advantage
melts away when increasing the overlap or the ansatz order.
While these properties are consistently observed with the GLL basis, the GL
results follow a less regular pattern and exhibit mostly lower convergence rates.
This behavior can partly be explained by the fact that, with Gauss points, using
a geometrically specified overlap width may result in zero overlapped node layers.
While this phenomenon appears only at low orders, the latter are always present
in the multigrid scheme, even at high ansatz orders.
A remedy to this problem is to apply a lower bound of ${\no[,l] = 1}$ for the nodal
overlap on every mesh level $l$.
Nevertheless, the GL-based approach remains slightly less efficient in comparison
with the GLL approach.
Therefore, further discussion will be constrained to the latter.

\begin{table}
\centering
\caption{Convergence rates obtained with different multigrid methods on isotopic
meshes composed of $(128/P)^3$ elements for increasing ansatz order $P$.}
\label{tab:isotropic:overlap}
\small
\begin{tabular}{r@{\hspace{1.25em}}l@{\hspace{1.25em}}l@{\hspace{1.25em}}l@{\hspace{1.25em}}c@{\hspace{1.25em}}c@{\hspace{.75em}}c@{\hspace{.75em}}c@{\hspace{.75em}}c}
\toprule
   &       &        &                                       &       &\multicolumn{4}{c}{$\bar r = -\lg\rho$} \\
\cline{6-9}
\# & Basis & Solver & Overlap                               & $\no$ & $P\!=\!4$ & $P\!=\!8$ & $P\!=\!16$ & $P\!=\!32$  \\
\midrule
 1 & GLL   & MG     & $\no = 1$                             &  $1,1,1,1, 1$ &  1.08  &  0.89  &  0.46  &  0.24  \\
 2 & GLL   & MG-CG  & $\no = 1$                             &  $1,1,1,1, 1$ &  1.34  &  1.21  &  0.80  &  0.52  \\
 3 & GLL   & MG     & $\deltaO = 0.08\Delta x$              &  $1,1,2,3, 6$ &  1.08  &  1.59  &  2.13  &  2.82  \\
 4 & GLL   & MG-CG  & $\deltaO = 0.08\Delta x$              &  $1,1,2,3, 6$ &  1.34  &  1.85  &  2.24  &  2.81  \\
 5 & GLL   & MG     & $\deltaO = 0.50\Delta x$              &  $2,3,5,9,17$ &  2.26  &  2.76  &  3.18  &  3.66  \\
 6 & GLL   & MG-CG  & $\deltaO = 0.50\Delta x$              &  $2,3,5,9,17$ &  2.29  &  2.67  &  3.08  &  3.75  \\
\midrule
 7 & GL    & MG     & $\no = 1$                             &  $1,1,1,1, 1$ &  0.73  &  0.96  &  0.73  &  0.40  \\
 8 & GL    & MG-CG  & $\no = 1$                             &  $1,1,1,1, 1$ &  1.49  &  1.26  &  1.07  &  0.78  \\
 9 & GL    & MG     & $\deltaO = 0.08\Delta x$              &  $0,1,1,3, 6$ &  0.73  &  1.15  &  1.70  &  2.20  \\
10 & GL    & MG-CG  & $\deltaO = 0.08\Delta x$              &  $0,1,1,3, 6$ &  1.14  &  1.38  &  1.84  &  2.15  \\
11 & GL    & MG     & $\deltaO = 0.50\Delta x$              &  $1,2,4,8,16$ &  1.83  &  2.28  &  1.52  &  0.94  \\
12 & GL    & MG-CG  & $\deltaO = 0.50\Delta x$              &  $1,2,4,8,16$ &  1.90  &  2.34  &  1.88  &  1.37  \\
13 & GL    & MG     & $\deltaO = 0.09\Delta x$, $\no \ge 1$ &  $1,1,2,3, 6$ &  1.36  &  1.75  &  1.87  &  2.21  \\
14 & GL    & MG-CG  & $\deltaO = 0.09\Delta x$, $\no \ge 1$ &  $1,1,2,3, 6$ &  1.49  &  1.81  &  2.03  &  2.22  \\
\bottomrule
\end{tabular}
%$^a$ Table foot note (with superscript)
\end{table}

\begin{figure}%[b]
\centering
\includegraphics[width=90mm]{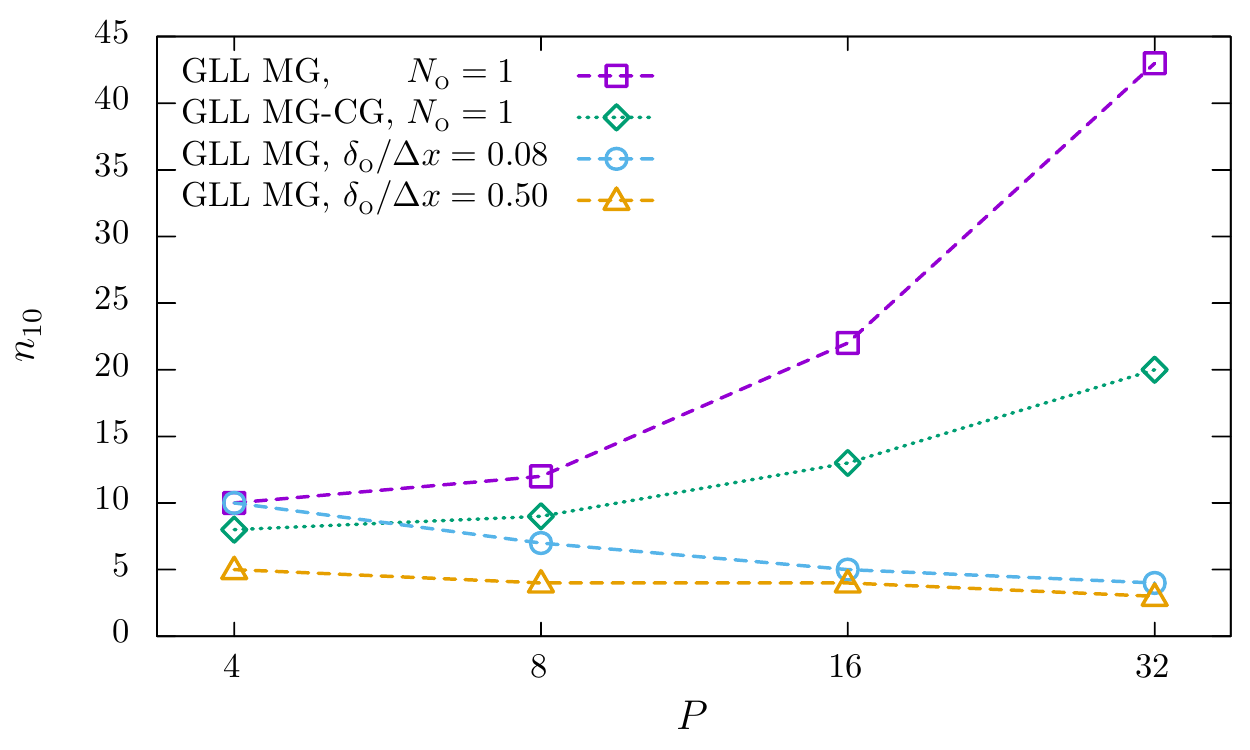}
\caption{Number of cycles required to reduce the residual by ten orders of magnitude for
selected methods listed in Tab.~\ref{tab:isotropic:overlap}.
}
\label{fig:p_n10}
\end{figure}

\begin{figure}%[b]
\centering
\includegraphics[width=90mm]{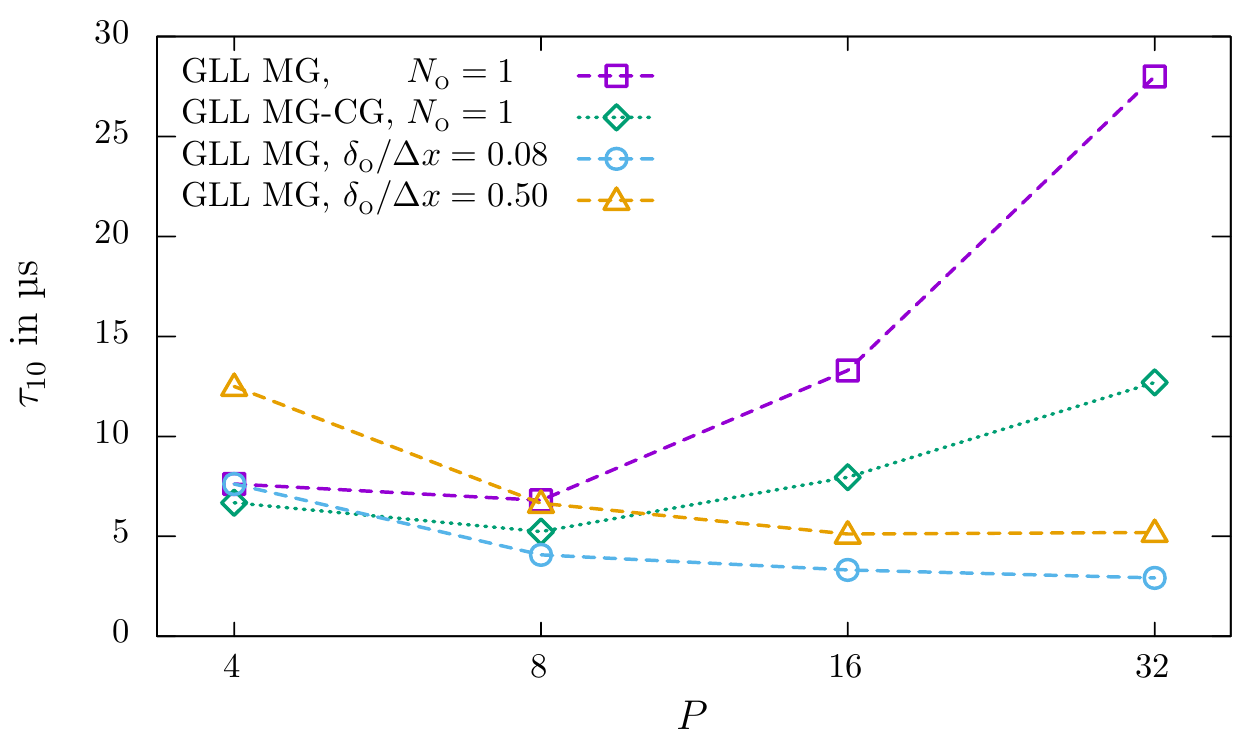}
\caption{Runtime per unknown required to reduce the residual by ten orders of magnitude for
selected methods listed in Tab.~\ref{tab:isotropic:overlap}.
}
\label{fig:p_tau10}
\end{figure}

\begin{figure}%[b]
\centering
\includegraphics[width=90mm]{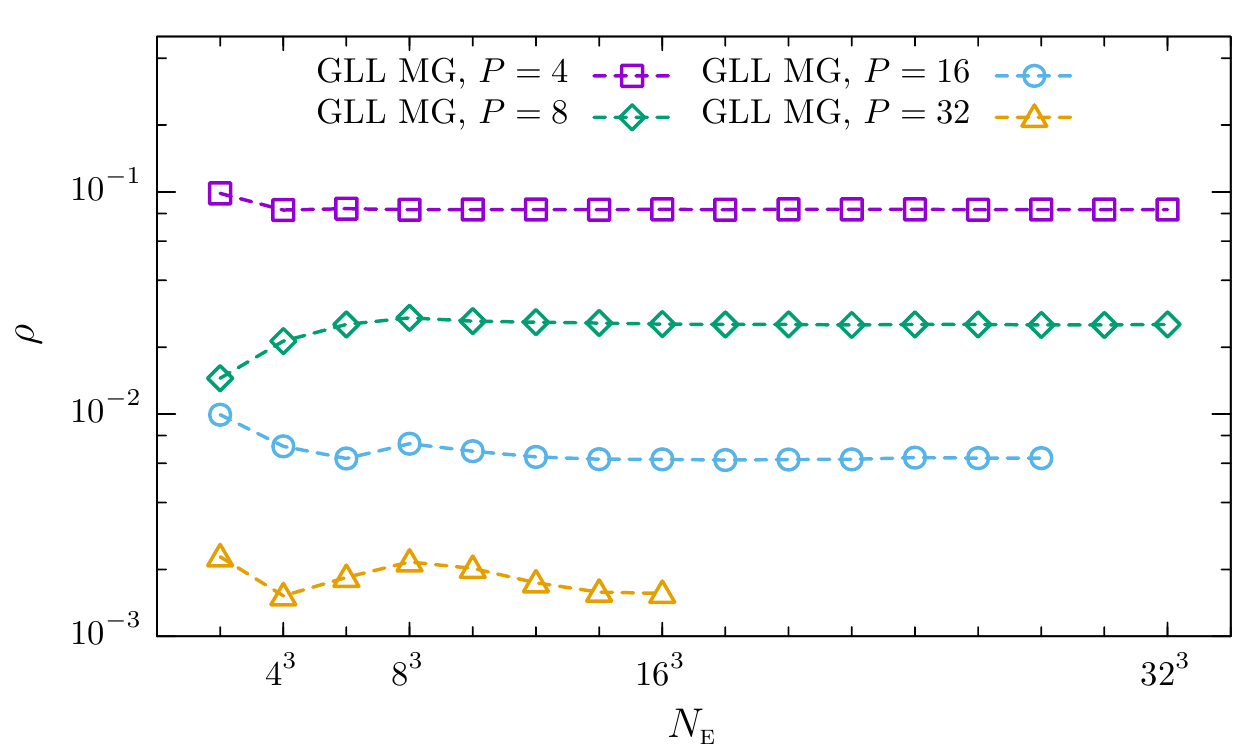}
\caption{Convergence rates for GLL MG with overlap ${\deltaO=0.08\Delta x}$
for different orders $P$ displayed as a function of the mesh size.
}
\label{fig:ne_rho}
\end{figure}

Complementary to the tabulated results, Fig.~\ref{fig:p_n10} depicts the number $n_{10}$ of multigrid cycles
that are required to reduce the Euclidian residual norm by ten orders of magnitude for selected cases
listed in Tab.~\ref{tab:isotropic:overlap}.
In agreement with the above discussion, the cycle count increases considerably when using GLL MG with
only one node layer overlap.
Adding Krylov acceleration (MG-CG) ameliorates this drawback, especially at higher order.
Yet, pure MG with a fixed geometric overlap of 8 percent is far more efficient and even
attains a decreasing $n_{10}$ with growing $P$.
As expected, using an overlap of $\Delta x$$/$$2$ yields a further reduction of the cycle count.
This advantage is, however, bought with additional computational cost related to the larger subdomain
operator $\NM A_{ss}$.
Figure~\ref{fig:p_tau10} confirms that GLL MG with ${\deltaO=0.08\Delta x}$ outpaces the other
choices for all polynomial degrees but 4, where the Krylov-accelerated method (MG-CG) with
one node overlap is slightly faster.
Moreover, Fig.~\ref{fig:ne_rho} illustrates the robustness of this method with respect to the
mesh size. With ansatz orders up to 16, the convergence rate becomes mesh independent
for ${\N{e} \gtrsim 12^3}$, whereas it still tends to improve beyond ${\N{e}=16^3}$ for ${P=32}$.
It is further worth noting that the convergence rates improves with growing order, reaching
an excellent ${\rho\approx6.3\times10^{-3}}$ with ${P=16}$ and
even better ${\rho\approx1.6\times10^{-3}}$ with ${P=32}$.
Moreover, runtimes of about $3.5$\,{\textmu}s per degree of freedom allow to solve problems
up to a million unknowns conveniently on a single core.

%-------------------------------------------------------------------------------

\subsection{Anisotropic Meshes}

\begin{figure}[p]

\centering

\subfloat[Multigrid convergence rates depending on the aspect ratio for
methods listed in Tab.~\ref{tab:ar:cases}.
\label{fig:ar_rho}]{\hspace*{4em}\includegraphics[width=90mm]{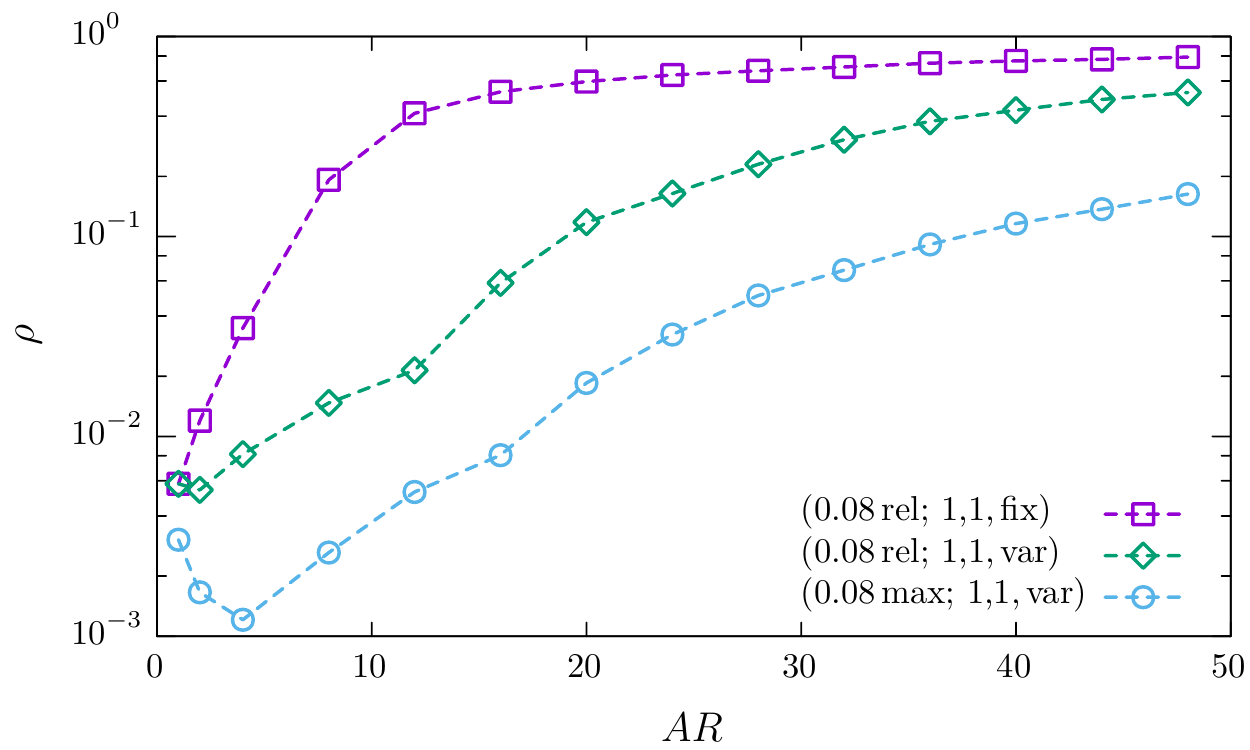}\hspace*{4em}}

%\medskip

\subfloat[Cycle counts depending on the aspect ratio for
methods listed in Tab.~\ref{tab:ar:cases}.
\label{fig:ar_n10}]{\hspace*{4em}\includegraphics[width=90mm]{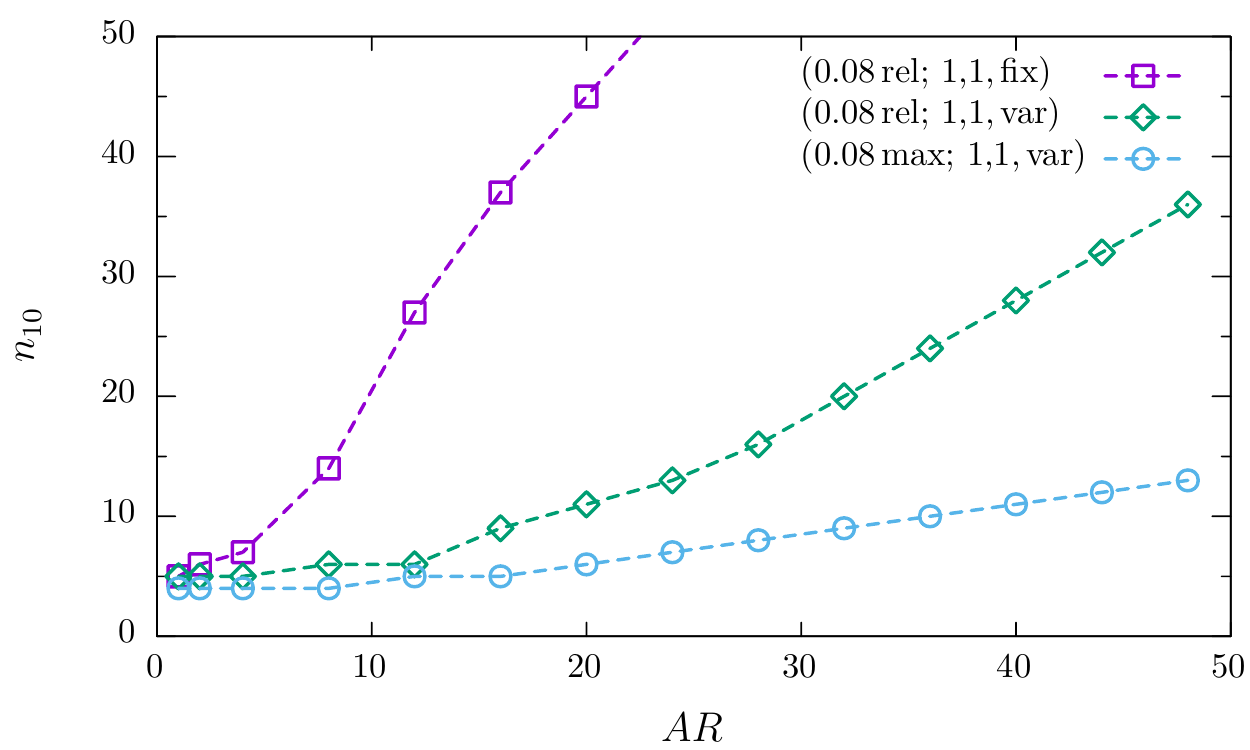}\hspace*{4em}}

%\medskip

\subfloat[Runtime per unknown depending on the aspect ratio for
methods listed in Tab.~\ref{tab:ar:cases}.
\label{fig:ar_tau10}]{\hspace*{4em}\includegraphics[width=90mm]{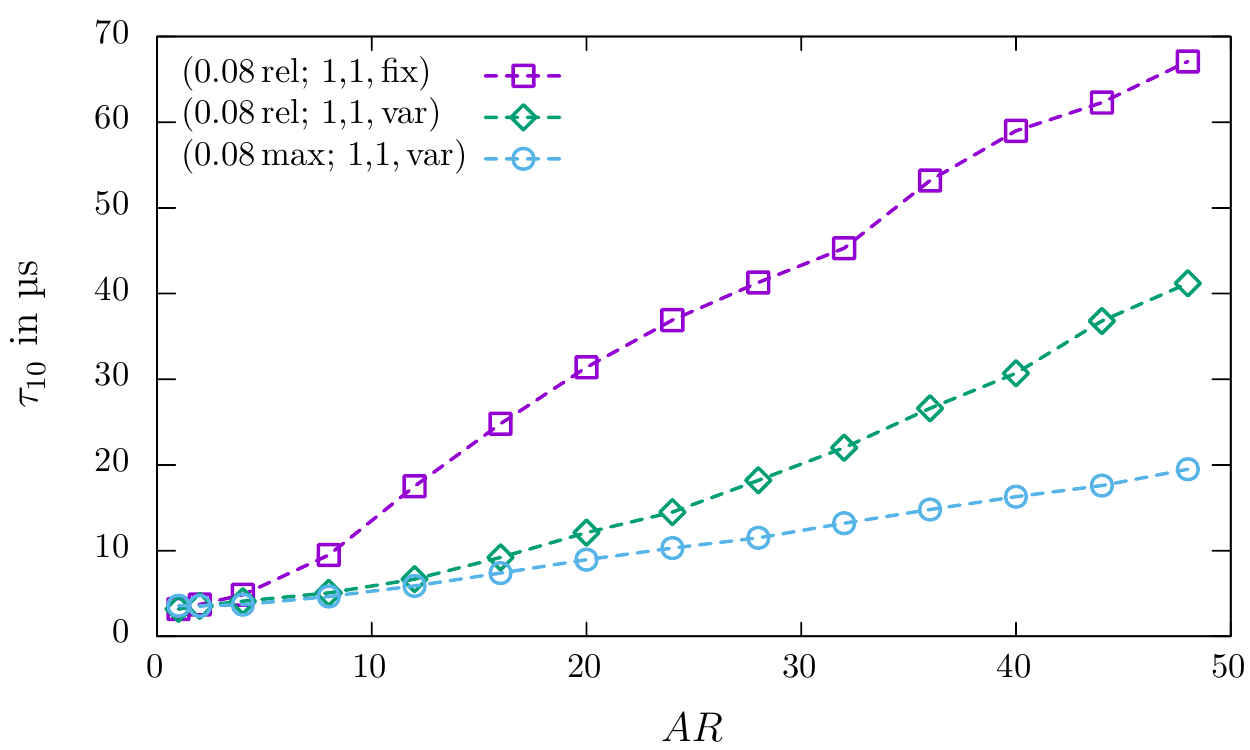}\hspace*{4em}}

\caption{Results for varying element aspect ratios.
\label{fig:ar}}

\end{figure}

As a second issue we investigated the suitability of the approach for anisotropic
meshes.
For this purpose, we defined a sequence of domains
\[
  \Omega = (0, 2 \pi\, AR) \times (0,2\pi\lceil\, AR/2 \rceil) \times (0,2\pi),
  \quad
\]
with aspect ratios $AR$ ranging from 1 to 48.
Using a uniform tessellation featuring the same number of elements in each coordinate
direction, $AR$ also represents the maximum element side aspect ratio.
Thus, for example, ${AR=32}$ results in
${\Delta x = 2\Delta y = 32\Delta z}$.
In earlier 2D studies, Krylov acceleration and variable V-cycles proved helpful,
though yet insufficient, for coping with anisotropy.
Based on this experience, we selected methods with different overlap and smoothing strategies, which are summarized in Tab.~\ref{tab:ar:cases}.
Methods
(0.08\,rel; 1,1,\,fix) and
(0.08\,rel; 1,1,\,var)
both use a relative subdomain overlap of 8 percent, which means that
the overlap width varies in each coordinate direction proportionally
to the corresponding element extension.
In contrast,
(0.08\,max; 1,1,\,var) sets the overlap width to 8 percent of the maximal
side length, but not larger than the element width in the given direction.
Additionally, the last two methods apply a variable V-cycle, which increases the
number of smoothing steps by a factor of 3 with each coarser level.
The performance of these methods was studied on a $8^3$ tessellation using
elements of order ${P=16}$.
Figure~\ref{fig:ar} depicts the obtained convergence rates, cycle counts and
runtimes per unknown for aspect ratios up to 48.
With method (0.08\,rel; 1,1,\,fix) convergence starts to degrade
at moderate aspect ratios and has already slowed by to orders of magnitude at
${AR=16}$.
Using a variable V-cycle improves the robustness such that a nearly constant
cycle count $n_{10}$ is maintained until ${AR=12}$.
From here convergence degrades more quickly, but remains superior to the
previous case.
Setting the overlap proportional to the largest element extension, as with
method (0.08\,max; 1,1,\,var), yields a further improvement,
which becomes even more pronounced in the range ${AR > 12}$, where the
overlap in the most compressed direction is already constrained by the
element width.
Compared to the isotropic case, the cycle count grew from 4 to 13 at
${AR=48}$, whereas the serial runtime increased by factor of 5.4 to
approximately 19\,{\textmu}s per unknown.
This seems to be a good starting point, given the prospect of further
acceleration, e.g. by parallelization.

\begin{table}
\renewcommand{\arraystretch}{1.5}
\caption{Test cases for investigating the robustness against the element aspect ratio.}
\label{tab:ar:cases}
\begin{tabular}{l@{\hspace{1.75em}}l@{\hspace{1.75em}}l}
\toprule
Case & Subdomain overlap & Smoothing steps
\\
\midrule
(0.08\,rel; 1,1,\,fix)
& $\deltaO[,\,x_i]=0.08\Delta x_i$
& $\ns[,l] = (1,1)$
\\
(0.08\,rel; 1,1,\,var)
& $\deltaO[,\,x_i]=0.08\Delta x_i$
& $\ns[,l] = (1,1) \times 3^{L-l}$
\\
(0.08\,max; 1,1,\,var)
& $\deltaO[,\,x_i] = \min[\max_j(0.08\Delta x_j),\,\Delta x_i]$
& $\ns[,l] = (1,1) \times 3^{L-l}$
\\
\bottomrule
\end{tabular}
\end{table}

%===============================================================================

\section{Conclusions}
\label{sec:concl}

We developed a polynomial multigrid method for nodal interior-penalty
formulations of the Poisson equation on three-dimensional Cartesian grids.
Its key ingredient is an overlapping weighted Schwarz smoother, which exploits
the underlying tensor-product structure for fast solution of the subdomain
problems.
The method achieves excellent convergence rates and proved robust against the
mesh size and ansatz orders up to at least 32.
Extending the ideas put forward in \cite{Stiller2016}, we showed that combining
Krylov acceleration, variable smoothing and increasing the subdomain overlap
proportionally to the maximum element width improves the robustness considerably
and renders the approach feasible for aspect ratios up to 50.
Moreover, the method is computationally efficient, allowing to solve problems
with a million unknowns in a few seconds on a single CPU core.

%===============================================================================

\section*{Acknowledgement}
Funding by German Research Foundation (DFG) in frame of the project
\mbox{STI~157/4-1} is gratefully acknowledged.

%%%%%%%%%%%%%%%%%%%%%%%%%%%%%%%%%%%%%%%%%%%%%%%%%%%%%%%%%%%%%%%%%%%%%%%%%%%%%%%%
% Bibliography

\end{document}